\magnification=\magstep1
\hsize16truecm
\vsize21.5truecm
\topskip=1truecm
\raggedbottom
\abovedisplayskip=3mm
\belowdisplayskip=3mm
\abovedisplayshortskip=0mm
\belowdisplayshortskip=2mm
\normalbaselineskip=12pt
\normalbaselines
\font\titlefont= cmcsc10 at 12pt
\def\F{\Bbb F}

\def\P{\Bbb P}

\input amssym.def
\input amssym.tex

\vskip 6.5pc
\noindent
\font\eighteenbf=cmbx10 scaled\magstep3
\font\titlefont=cmcsc10 at 12pt
\vskip 2.0pc
\centerline{\eighteenbf  Generalized Reed-Muller Codes and}
\vskip 1.5pc
\centerline{\eighteenbf Curves with Many Points}
\vskip 3.0pc
\centerline {\titlefont G. van der Geer and M. van der
Vlugt}
\vskip 2.0pc
\noindent 
In the quest for curves over finite fields with many points coding theory
has been a useful guide, as words of low weight in trace codes correspond to
Artin-Schreier curves with many points. This correspondence can be
extended to subcodes of low weight and fibre products of Artin-Schreier
curves. Subcodes of minimum weight of a code
${\cal C}$ determine the weight hierarchy of ${\cal C}$ and knowledge of the
weight hierarchy indicates where curves with many points are likely to be
found. However, determination of weight hierarchies is a hard problem in
coding theory. In 1990 Wei found the weight hierarchy of the classical
binary Reed--Muller codes (see [W]). Six years later Heijnen and Pellikaan
succeeded in finding the weight hierarchy of Reed--Muller codes over
arbitrary finite fields, cf.\ [H-P]. In this paper we derive some results on
curves with many points which are closely related to the weight hierarchy of
the $q$-ary or generalized Reed-Muller codes. We use subcodes of
generalized Reed-Muller codes to construct our curves. We also present some
types of curves which attain the Hasse-Weil upper bound.
\vskip\baselineskip
{\bf 1. On the weight hierarchy of generalized Reed-Muller codes}
\vskip\baselineskip
\noindent
Let ${\cal C}$ be a (linear) code of length $n$ and dimension $k$ over a
finite field
$\F_{q=p^m}$. An important parameter of a subcode ${\cal D}$ is its weight
$w({\cal D})$, by which we mean the number of coordinate places for which at
least one word of ${\cal D}$ has a non-zero coordinate. Since the projection
of ${\cal D}$ onto a coordinate place is a $\F_q$-linear map we have
$$
w({\cal D}) = {1\over q^r-q^{r-1}}\sum_{d\in {\cal D}} w(d), \eqno(1)
$$
where $r=\dim({\cal D})$ and $w(d)$ is the weight of the word $d$. The $r$-th
generalized Hamming weight $d_r({\cal C})$ for $1 \leq r \leq k$ is defined
by
$$
d_r({\cal C})=\min\{ w({\cal D}): {\cal D} \hbox{~is an $r$-dimensional
subcode of ${\cal C}$}\}.
$$
The set $\{d_r({\cal C}) : 1 \leq r \leq k\}$ is called the {\sl weight
hierarchy} of ${\cal C}$.
\vskip\baselineskip
The  family of $q$-ary or generalized Reed-Muller
codes can be defined as follows. Elements of the vector space
$$
P_s = \{f\in \F_q[X_1,\ldots,X_m]:\deg (f) \leq s\}
$$
can be evaluated at the points of the affine space $\F_q^m$. This defines
an evaluation map
$$\beta : P_s \rightarrow \F_q^n
$$ 
with $n=q^m$ given by $f \mapsto(f(v)_{v\in\F_q^m})$.   Note that the kernel
of $\beta$ is the ideal generated by the polynomials $X_i^q-X_i$.
The image $\beta(P_s)$ is the $q$-ary Reed-Muller code $R_q(s,m)$ of order
$s$ in $m$ variables. One thus obtains a large class of codes. 
\par
\smallskip
In the paper by Heijnen and Pellikaan [H-P] one can find the following
algorithm which  determines an $r$-dimensional subcode of
$R_q(s,m)$ of minimum weight. First we fix an enumeration $\{\alpha_0,
\alpha_1,\ldots,\alpha_{q-1}\}$ of $\F_q$.
\parindent14pt
\item{1.} Let $Q= \{ 0,1,\ldots,q-1\}$. The set $Q^m$ is ordered
lexicographically using the natural order on $Q$. Elements of $Q^m$ are
denoted by $\sigma=(i_1,\ldots,i_m)$ and have a degree
defined by
$$
\deg (\sigma)=\sum_{j=1}^m i_j .
$$
\item{}
\item{2.} To an element  $\sigma\in Q^m$  associate  the following 
polynomial $f=f_{\sigma}$ (which is {\sl reduced} modulo the kernel of
$\beta$): 
$$
f = \prod_{j=1}^m \prod_{t=i_j+1}^{q-1} (X_j -
\alpha_t).\eqno{(2)}
$$
\item{3.} Take the first $r$ elements of degree $\geq m(q-1)-s$ in $Q^m$.
\par
\smallskip
\noindent According to one of the main results 
in [H-P] (Thm.5.9) the outcome of the algorithm is:
\smallskip
\noindent
\proclaim (1.1) Theorem.  Let $\sigma_1, \ldots,\sigma_r$ be the first
$r$ elements in $Q^m$ with $\deg(\sigma_i) \geq m (q-1)-s$. Then the
codewords induced by the polynomials $f_1,\ldots, f_r$ associated to
$\sigma_1,\ldots, \sigma_r$ generate an $r$-dimensional subcode of minimum
weight of $R_q (s,m)$.\par
\medskip
The same theorem in [H-P] also yields a formula for $d_r(R_q(s,m))$.
\vskip\baselineskip
\proclaim (1.2) Formula.  If $\sigma_r=(i_1,\ldots,i_m)$ is the $r$-th element
in
$Q^m$ of degree $\geq m(q-1)-s$  then
$$
d_r(R_q(s,m)) = 1 + \sum^m_{j=1} i_{m-j+1} q^{j-1} .
$$
\par
\smallskip
To illustrate this algorithm we consider an example.
\smallskip\noindent
{\bf (1.3) Example.} We enumerate $\F_ p = \{\alpha_0, \alpha_1,
\ldots, \alpha_{p-1}\} = \{p-1, p-2, \ldots, 0\}$ and we consider $R_p(2,3)$
for an odd prime $p$. The first four elements $\sigma \in Q^m$ with degree
$\geq 3(p-1)-2 = 3p-5$ are $(p-3,p-1, p-1)$, $(p-2, p-2, p-1)$, 
$(p-2, p-1, p-2)$ and $(p-2, p-1, p-1)$.
\par
According to (2) the corresponding $f_i$ are
$$
f_1 = (X_1-1)X_1,\quad f_2 = X_1X_2,\quad f_3 = X_1X_3 \, {\rm~and~} \,  f_4
= X_1 .
$$
Now the codewords induced by $f_1 = X^2_1 - X_1$, $f_2 = X_1X_2$ and $f_3
= X_1X_3$ generate a 3-dimensional subcode of $R_p(2,3)$ denoted by $\langle 
f_1, f_2, f_3\rangle $ and 
$$
d_3(R_p(2,3)) = 1 + (p-2) + (p-1) p + (p-2)p^2 =
p^3-p^2-1.
$$
\smallskip
In the next section we relate curves to codewords in $R_q(s,m)$ and
then we apply the results of this section to obtain results for curves.
\bigskip
\centerline{\bf 2. Words in $R_q(s,m)$ and Artin-Schreier curves}
\medskip
\noindent
>From now on we identify reduced polynomials $f$ in $P_s$ and the
codewords $c_f$ which they induce by the evaluation map. A polynomial
$f$ of the form (2) is a product of independent (inhomogenenous) linear
forms $(X_j - \alpha_t)$. For fixed $a \in \F_{q^m}$ we consider
${\rm Tr}(ax)$ where ${\rm Tr}= {\rm Tr}_{q^m/q}$ is the trace map from
$\F_{q^m}$ onto
$\F_q$. When we view
$\F_{q^m}$ as an $m$-dimensional vector space over $\F_q$ then ${\rm
Tr}(ax)$ is a linear form in $\F_q[x_1, \ldots, x_m]$. If we choose
$\F_q$-independent $a_1, a_2,\ldots, a_m$ in $\F_{q^m}$ we can write (2)
as
$$
f = \prod_{j=1}^m \prod_{t=i_j+1}^{q-1} ({\rm Tr}(a_jx)-\alpha_t) .
$$
Observe that for $a,b \in \F_{q^m}$ we have
$$
{\rm Tr}(ax) {\rm Tr}(bx) = {\rm Tr}({\rm Tr}(ax)bx) = {\rm
Tr}(\sum^{m-1}_{j=0} a^{2^j} b x^{2^j+1}).\eqno{(3)}
$$
Repeating this process we see that the codeword corresponding to $f$ can
be written in trace form:
$$
c_f = ({\rm Tr}(R(x)))_{x \in \F_{q^m}}{\rm~with~} R(x) \in \F_{q^m}[x].
$$
An element of $\F_{q^m}$ has trace zero precisely if it is of the form
$y^q-y$ for some $y \in \F_{q^m}$. Now we associate to the codeword $c_f$ the
irreducible complete smooth curve $C_f$ over $\F_{q^m}$ given by the affine
equation
$$
y^q-y = R(x). \eqno{(4)}
$$
In the sequel we always assume that $\deg(R)$ is prime to $p$, the
characteristic of $\F_{q^m}$. Then the irreducible smooth projective
Artin-Schreier curve $C_f$ given by (4) has genus
$$
g(C_f) = (q-1)(\deg(R) - 1)/2 .\eqno{(5)}
$$
We  see immediately that there is a relation between the weight $w(c_f)$
of $c_f$ and the number of $\F_{q^m}$-rational points on $C_f$:
$$
w(c_f) = q^m - (\# C_f(\F_{q^m}) - 1)/q . \eqno{(6)}
$$
This correspondence between codewords and curves can be extended to
subcodes and fibre product of curves (see  [G-V1]). If ${\cal D}$ is an
$r$-dimensional subcode of $R_q(s,m)$ with basis $c_{f_1}, \ldots,
c_{f_r}$ we consider the corresponding curves $C_{f_i}$ defined by $y^q - y
= R_i(x)$. Each of these curves admits a natural map $\varphi_i : C_{f_i}
\mapsto \P^1$. Then we associate to ${\cal D}$ the curve
$$
C_{\cal D} = {\rm Normalization~of~} C_{f_i}\times_{\P^1} \times \ldots
\times_{\P^1} C_{f_r}.
$$
>From [G-V1] we recall the following proposition.
\noindent
\proclaim
(2.1) Proposition. The weight $w({\cal D})$ satisfies $w({\cal D}) =
q^m - (\#C_{\cal D}(\F_{q^m}) - 1)/q^r$. \par
\smallskip
For $f$ in the $\F_q$-vector space generated by $f_1, \ldots, f_r$
we denote the trace of Frobenius on $C_f$ by $\tau_f$ (i.e.
$\# C_f (\F_{q^m}) = q^m + 1 - \tau_f)$.\par
According to [G-V 1] the trace of Frobenius  $\tau_{\cal D}$ of $C_{\cal D}$
and the genus $g(C_{\cal D})$ satisfy
$$
(q-1)\tau_{\cal D} = \sum_{f\in \langle f_1, \ldots, f_r \rangle  -
\{0\}}\tau_f\eqno(7)
$$
$$
(q-1)g(C_{\cal D}) = \sum_{f\in \langle  f_1, \ldots, f_r \rangle - \{0\}}
g(C_f).
\eqno(8)
$$
>From (6) and Proposition (2.1) it is clear that words and subcodes of
low weight correspond to curves with many rational points.
\par
In the next section we illustrate these ideas by some examples.
\vskip\baselineskip
\centerline {\bf  3. The code $R_p(2,3)$}
\bigskip
We consider $R_p(2,3)$ for an odd prime $p$. As we saw at the end of
Section 1 the first four polynomials of the form (2) are
$$
f_1 = (X_1-1)X_1,\quad f_2 = X_1X_2,\quad f_3 = X_1X_3 \, {\rm~and~}\,  f_4 =
X_1.
$$
The number of zeros in the codeword $c_{f_1}$ is $2p^2$, so its weight
satisfies
$w(c_{f_1}) = p^3 - 2p^2$. We can write $f_1 = ({\rm Tr}(x)-1){\rm Tr}(x)$,
where ${\rm Tr} = {\rm Tr}_{{p^3}/p}$.
Applying (3) we find
$$
f_1 ={\rm Tr}(x^2 +x^{p+1} + x^{p^2+1})-{\rm Tr}(x).
$$
Since the trace map on $\F_{p^3}$ satisfies
$$
{\rm Tr}(x^{p^2+1}) = {\rm Tr}(x^{p+1})\eqno{(9)}
$$
we find for the corresponding codeword
$$
c_{f_1}={\rm Tr}(2x^{p+1} + x^2 - x)_{x\in\F_{p^3}}.
$$
The curve $C_{f_1}$ which we now associate to $c_{f_1}$ is given by the
affine equation
$$
y^p-y = 2x^{p+1} + x^2-x.
$$
>From the number of zeros in $c_{f_1}$ we immediately obtain
$\#C_{f_1}(\F_{p^3}) = 2p^3 + 1$,
while according to (5) the genus satisfies $g(C_{f_1}) = (p-1)p/2$.
\bigskip
\noindent
\proclaim  (3.1) Result. For $p=3$ we obtain a curve over $\F_{27}$ of
genus 3 with 55 points. 
\par
This is quite close to the upper bound $58$.
\bigskip
\noindent
{\bf (3.2) Remark.} Note that by replacing monomials by lower
degree monomials (using Frobenius as in (9)) or by neglecting the monomials or
polynomials which give zero under the trace map we can reduce  the genus of the
curve which we associate to the codeword. This possibility of {\sl genus
reduction} makes curves much more attractive as building blocks for fibre
product curves with many points.
\bigskip
If we change to $f_2 ={\rm Tr}(x){\rm Tr}(ax)$ with $a \in \F_{p^3}-\F_p$
then along the same lines we find
$$
c_{f_2} ={\rm Tr}((a^p+a)x^{p+1}+ax^2)
$$
with $a^p+a\not= 0$. The word $c_{f_2}$ has $2p^2-p$ zeros and the corresponding
curve $C_{f_2}$ given by
$$
y^p - y = (a^p + a)x^{p+1} + ax^2
$$
has genus $(p-1)p/2$ and $\# C_{f_2}(p^3) = 2p^3 - p^2+1$.
Now we consider the 3-dimensional subcode
$$
{\cal D} = \langle  ({\rm Tr}(x) - 1){\rm Tr}(x),~{\rm Tr}(x){\rm Tr}(ax),~{\rm
Tr}(x){\rm Tr}(bx) \rangle
$$
with $\{1,a,b\} \subset \F_{p^3}$ independent over $\F_p$.
\smallskip
\noindent
\proclaim (3.3) Proposition. The curve $C_{\cal D}$ defined over $\F_{p^3}$
corresponding to the subcode ${\cal D}= \langle c_{f_1}, c_{f_2},
c_{f_3}\rangle$ has genus
$(p^4-p)/2$ and $\# C_{\cal D}(\F_{p^3}) = p^5 + p^3 + 1$.\par
\medskip
\noindent
{\sl Proof.} The curves which occur in the fibre product induced by ${\cal D}$
are of the form
$$
y^p-y = (2\lambda_1 + \lambda_2(a^p+a) + \lambda_3(b^p+b))x^{p+1} + (\lambda_1
+ \lambda_2 a + \lambda_3b)x^2 - \lambda_1x\eqno{(10)}
$$
with $(\lambda_1, \lambda_2, \lambda_3) \in \F^3_p - \{0\}$. The coefficient
of $x^{p+1}$ can be written as
$$
(\lambda_1 + \lambda_2a + \lambda_3b)^p + (\lambda_1 + \lambda_2a +
\lambda_3b).
$$
If this coefficient is zero then $\lambda_1 + \lambda_2a +
\lambda_3b
\in \F_{p^2} \cap \F_{p^3} = \F_p$. Since $\{1,a,b\}$ are
$\F_p$-independent we find $\lambda_2 = \lambda_3=0$ which implies $\lambda_1
= 0$. So the right-hand side of (10) has degree $p+1$ and the curves
defined by (10) have genus $(p-1) p/2$. From (7) we derive
$$
g(C_{\cal D}) = (p^3-1)p/2 .
$$
Moreover we saw at the end of Section 1 that $d_3(R_p(2,3)) = p^3-p^2-1$ and
then Proposition (2.1) yields $\# C_{\cal D}(\F_{p^3}) = p^5 + p^3 + 1$.
\smallskip
\noindent
\proclaim (3.4) Corollary. For $p=3$ the curve $C_{\cal D}$, which is defined
over $\F_{27}$, has genus 39 and $\# C_{\cal D} (\F_{27}) = 271$. \par
\vskip\baselineskip
Corollary (3.4) gives an improvement of the best value known until now for
$(q,g)=(27,39)$ which is 244 (see [G-V3]).
\smallskip
When we take the subcode ${\cal D} = \langle  c_{f_1}, c_{f_2}, c_{f_3}, c_{f_4}
\rangle$ we add the word corresponding to $f_4 = X_1$ to the basis of the
former subcode. Note that the associated curve $C_{f_4}$ given by $y^p-y=x$
has genus zero. Following the above procedure we find:
\smallskip
\noindent
\proclaim (3.5) Proposition. The curve $C_{\cal D}$ defined over $\F_{p^3}$
which corresponds to the subcode ${\cal D} = \langle c_{f_1}, c_{f_2}, c_{f_3},
c_{f_4}\rangle $ has genus $(p^4 - p)p/2$ and $\# C_{\cal D}(\F_{p^3}) =
p^6+1$. \par
\smallskip
\noindent
\proclaim (3.6) Corollary. For $p=3$ the curve $C_{\cal D}$ over $\F_{27}$ has
genus 117 and $\# C_{\cal D}(\F_{27}) = 730$.\par
\smallskip
This is fairly good compared to Oesterl\'e's upper bound which is 859
for curves of genus 117 over $\F_{27}$.
\vskip\baselineskip
The next example shows that we can also use subcodes which are not of minimum
weight to construct curves with many points. This is caused by the fact that
the possibility of genus reduction is a very useful feature.
\bigskip
\centerline{\bf 4. The code $R_3(3,3)$}
\bigskip
\noindent
Here we exploit the code $R_3(3,3)$ to produce a good curve 
over $\F_{27}$.  In this case the first three elements
$\sigma \in Q^3$ with $\deg(\sigma) \geq 3$ are $(0,1,2)$, $(0,2,1)$ and
$(0,2,2)$ with corresponding $f_1=(X_1-1)X_1X_2$, $f_2 =
(X_1-1)X_1X_3$ and $f_3 = (X_1-1)X_1$. We can write $f_1$ as
$$
f_1 = ({\rm Tr}(x)-1){\rm Tr}(x){\rm Tr}(ax) {\rm ~with~} \F_3{\rm
-independent}~ \{1,a\}
\subset \F_{27}
$$
and we find for the codeword
$$
c_{f_1}={\rm Tr}[ 2a^9x^{13} + (2a^3 + a)x^7 + (a^3 + 2a)x^5 - (a^3+a)x^4 +
ax^3 -ax^2 ].
$$
The word $c_{f_1}$ has 21 zeros and the related curve has genus 12, but if we
could eliminate $2a^9x^{13}$ we get a curve of genus 6 since $2a^3 + a\not= 0$.
As $x^{13}\in\F_3$ for $x\in\F_{27}$ we can neglect $2a^9x^{13}$ if we
require ${\rm Tr}_{{27}/3}(a)=0$. Unfortunately there is no
$\F_3$-independent subspace $\langle1,a,b\rangle $ in $\F_{27}$ with ${\rm
Tr}(a)={\rm Tr}(b)=0$. However, when we take the fibre product of the curves
induced by
$f_1$ and
$f_3$ we find:
\smallskip
\noindent
\proclaim (4.1) Proposition. For $a \in \F_{27} - \F_3$ with ${\rm
Tr}(a)=0$ the fibre product of the curves $C_{f_1}$ and $C_{f_3}$ has genus
$g=21$ and possesses 163 points over $\F_{27}$. \par
\smallskip
\noindent
{\sl Proof. } The curve $C_{f_1}$ is given by
$$
y^3-y=(2a^3+a)x^7+(a^3+2a)x^5-(a^3+a)x^4+ax^3-ax^2,
$$
while  $C_{f_2}$ is defined by
$$
y^3-y = 2x^4+x^2-x .
$$
The words in the subcode generated by $c_{f_1}$ and $c_{f_3}$ are
$$
\lambda_1c_{f_1} + \lambda_2 c_{f_3} = (X_1 - 1)X_1(\lambda_1X_2 +
\lambda_2)
$$ 
with $\lambda_1, \lambda_2 \in \F_3$. For $\lambda_1\not= 0$ we
have 21 zeros or weight 6 and for $\lambda_1=0$ the corresponding word has
weight 9. It follows from (1) that
$$
w({\cal D} = \langle  c_{f_1}, c_{f_3}\rangle) = (6\cdot 6 + 2\cdot 9)/6=9 .
$$
Then Proposition (2.1) implies $\# C_{\cal D}(\F_{27}) = 163$. The curves which
occur in the fibre product are of the form
$$
y^3-y = \lambda_1(2a^3+a)x^7+\lambda_1(a^3+2a)x^5 +
(2\lambda_2-\lambda_1(a^3+a))x^4 + \lambda_1ax^3 +
(\lambda_2-\lambda_1a)x^2-\lambda_2x
$$
with $(\lambda_1,\lambda_2) \in \F^2_3 - \{0\}$. If $\lambda_1 \not= 0$ the
genus is $6$, while for $\lambda_1 = 0$ the genus is $3$ and we deduce from 
(8) that $g(C_{\cal D}) = (36+6)/2=21$. $\quad \square$
\smallskip
The curve from Proposition (4.1) satisfies the conditions of entry to the
Tables [G-V3] since Ihara's upper bound in this case is 214.
\bigskip
\centerline{\bf 5. Curves that attain the Hasse-Weil upper bound}
\bigskip
\noindent
In this section we construct maximal curves, i.e.\ irreducible smooth
curves of genus $g$ over $\F_q$ for which the number of $\F_q$-rational
points  attains the Hasse-Weil upper bound $q+1+2g\sqrt{q}$. This can only
happen if the genus is small compared to $q$, more precisely, by [S-X] and
[F-T] we know that for such curves $g\leq (\sqrt{q}-1)^2/4$
or $g=(q-\sqrt{q})/2$. The curves we construct here arise from the Reed-Muller
codes ${R}_p(2,m)$ for odd primes $p$ and even number of variables $m$.

According to the algorithm in Section 1 the polynomials (2) corresponding to
the first $m$ elements in $Q^m$ with degree $\geq m(p-1)-2$ are
$$
f_1=(X_1-1)X_1=({\rm Tr}(x)-1){\rm Tr}(x),\quad f_i=X_1X_i={\rm Tr}(x){\rm
Tr}(a_ix)\qquad (2\leq i \leq m),
$$
where the $a_i$ are chosen such that $\{ 1,a_2,\ldots,a_m\}$ are
$\F_p$-independent elements of $\F_{p^m}$. To these polynomials $f_i$ we can
associate Artin-Schreier curves. For $f_1$ we find
$$
C_{f_1}: \qquad y^p-y=x^{p^{m/2}+1}+2(\sum_{j=1}^{(m/2)-1} x^{p^j+1}) +
x^2-x,
$$
a curve with genus $g(C_{f_1})=p^{m/2}(p-1)/2$ and $\#
C_{f_1}(\F_{p^m})=2p^m+1$. The elements $f_i$ with $2\leq i \leq m$ define
curves
$$
C_{f_i}: \qquad
y^p-y=a_i^{p^{m/2}}x^{p^{m/2}+1}+\sum_{j=1}^{(m/2)-1}(a_i^{p^j}+a_i)
x^{p^j+1} + a_ix^2. \eqno(11)
$$
Now we observe that for all $x \in \F_{p^m}$ we have
$$
{\rm Tr}_{p^m/p}(a_i^{p^{m/2}}x^{p^{m/2}+1})={\rm
Tr}_{p^{m/2}/p}(a_i^{p^{m/2}}+a_i)x^{p^{m/2}+1},
$$
so that if $a_i^{p^{m/2}}+a_i={\rm Tr}_{p^m/p^{(m/2)}}(a_i)=0$ 
the curve in $(11)$
has the same number of points as the Artin-Schreier curve with the term
$a_i^{p^{m/2}}x^{p^{m/2}+1}$ deleted. In this way we can reduce the genus
without changing the number of points.
\bigskip
\noindent
\proclaim (5.1) Proposition. For $q=p^m$ with $p$ odd, $m$ even and 
$1\leq r \leq m/2$ there exists a maximal curve $C$ over $\F_q$ with
$$
g(C)=p^{(m/2)-1}(p^r-1)/2 \qquad {\rm and}
\quad \#C(\F_q)=p^m+1+(p^r-1)p^{m-1}.
$$
\par
\smallskip
\noindent
{\sl Proof.} Using the observation just made we consider the subspace $L$ of
elements $a_i\in \F_{p^m}$ with ${\rm Tr}_{p^m/p^{(m/2)}}(a_i)=0$. This is a
subspace of dimension
$m/2$ not containing $1$. Now we take a basis $\{a_2,a_3,\ldots,a_{(m/2)+1}\}$
of $L$ and we consider the $r$-dimensional subcode $ {\cal D}=\langle
c_{f_2},\ldots,c_{f_{r+1}}\rangle $ of $R_p(2,m)$ with $1\leq r \leq m/2$. The
curves involved in the fibre product defined by ${\cal D}$ are
$$
y^p-y=\sum_{j=1}^{(m/2)-1}(\sum_{i=2}^{r+1}\lambda_i(a_i^{p^j}+a_i))x^{p^j+1}
+
(\sum_{i=2}^{r+1}\lambda_ia_i)x^2\eqno(12)
$$
with $(\lambda_2,\ldots,\lambda_{r+1})\in \F_p^r-\{0\}$. Since the elements 
$a_2,\ldots,a_{r+1}$ are independent over $\F_p$ and satisfy
$a_i^{p^{m/2}}+a_i=0$ the coefficient of $x^{p^{(m/2)-1}+1}$ in (12) is not
zero. This implies that the genus of the curves given by (12) is
$(p-1)p^{(m/2)-1}/2$. It follows from (8) that
$$
g(C_{\cal D})=p^{(m/2)-1}(p^r-1)/2.
$$
Furthermore, as the curves (12) are related to 
$X_1(\lambda_2X_2+\ldots+\lambda_{r+1}X_{r+1})$, which  has $2p^{m-1}-p^{m-2}$
zeros, we find
$$
\# C_f(\F_{p^m})=p^m+1+(p^m-p^{m-1}).
$$
By formula (7) we obtain $\#C_{\cal D}(\F_q)=p^m+1+(p^r-1)p^{m-1}$.
\smallskip
\noindent
{\bf (5.2) Remark.} From the proof of Proposition (5.1) we deduce three other
families of maximal curves.
\item{$\bullet$} Curves of the form
$$
y^p-y=(\sum_{j=1}^{(m/2)-1}(a^{p^j}+a)x^{p^j+1})+ax^2\qquad {\rm with}\quad
{\rm Tr}_{p^m/p^{(m/2)}}(a)=0
$$
are maximal curves over $\F_{p^m}$ of genus $p^{(m/2)-1}(p-1)/2$.
\par
\item{$\bullet$}  Curves of the form
$$
y^p-y=a^{p^{m/2}}x^{p^{(m/2) }+1} \qquad {\rm with}\quad
{\rm Tr}_{p^m/p^{(m/2)}}(a)=0
$$
have $p^{m+1}+1$ rational points over $\F_{p^m}$ and genus $p^{m/2}(p-1)/2$,
so these are maximal too.
\par
\item{$\bullet$} For $1\leq r \leq m/2$ an $r$-dimensional fibre product of
curves from the preceding type  yields a maximal curve of genus
$(p^r-1)p^{m/2}/2$ with $p^{r+m}+1$ points.
\par
\noindent
These maximal curves were mentioned in [G-V2].
\bigskip
As was remarked in [L], any quotient $C^{\prime}$ of a maximal curve $C$ over
$\F_q$ which is defined over $\F_q$ (i.e.\ the map $C\to C^{\prime}$ is
defined over $\F_q$) is automatically maximal. This enables one to obtain
new maximal curves from the above ones.
\bigskip
\noindent
\proclaim (5.3) Proposition. For $q=p^m$ with $m$ even and $r$ a divisor of
$m$ with $1 \leq r \leq m/2$ there exists a maximal curve $C^{\prime}$ over
$\F_q$ with
$$
\eqalign{ g(C^{\prime })&=(p^{(m/2)-1}-1)(p^r-1)/4,\cr
\#C^{\prime}(\F_q)&=p^m+1+(p^{m-1}-p^{m/2})(p^r-1)/2.\cr}
$$
\par
\smallskip
\noindent
{\sl Proof.} Consider the maximal curve $C$ from Proposition (4.1) which is the
fibre product of the Artin-Schreier curves
$$
y_i^p-y_i=\sum_{j=1}^{(m/2)-1}(a_i^{p^j}+a_i)x^{p^j+1}+a_ix^2\qquad {\rm
for}\quad 2\leq i \leq r+1.
$$
Since each of these admits the involution
$(x,y_i) \mapsto(-x,y_i)$ we have an involution on $C$, the fixed points of
which are the $p^r$ points lying over
$x=0$ and the point over $x=\infty$. Applying the Hurwitz-Zeuthen formula 
$$
2g(C)-2= 2(2g(C^{\prime})-2)+p^r+1
$$
we find the
genus of the quotient curve 
$$
g(C^{\prime})= (p^{(m/2)-1}-1)(p^r-1)/4.
$$
The maximality implies that the number of points is as given. $\square$
\smallskip
Maximal curves of the form
$$
y^p-y=a^{p^{m/2}}x^{p^{(m/2)+1}}\qquad {\rm with} \quad {\rm
Tr}_{p^m/p^{m/2}}(a)=0
$$
possess an action  of the $d$-th roots of unity defined  via
$(x,y) \mapsto (\zeta x, y)$ provided $d$ divides $p^{m/2}+1$. The fixed
points under these automorphisms are again the points lying over $x=0$ and
the point at infinity. The ramification is tame of order $d$ and this gives:
\bigskip
\noindent
\proclaim (5.4) Proposition. For any positive divisor $d$ of $p^{m/2}+1$
there exists a maximal curve $C^{\prime}$ over $\F_{p^m}$ with
$$
\eqalign{ g(C^{\prime})& =(p^{m/2}-d+1)(p-1)/2d, \cr
\# C^{\prime}(\F_{p^m})& =p^m+1+(p^m-(d-1)p^{m/2})(p-1)/d.\cr}
$$
\par
\smallskip
If we take a fibre product $C$ as sketched in the third point of Remark (5.2)
with $r$ a divisor of $m$ and $1\leq r \leq m/2$ we find that $C$ is given by
an equation
$$
y^{p^r}-y=F(x^{p^{(m/2)}+1}) \qquad {\rm with}\qquad F(t) \in \F_{p^m}[t].
$$
The $d$-th roots of unity act as remarked above by which we obtain the following
result.
\bigskip
\noindent
\proclaim (5.5) Proposition. For any positive divisor $d$ of $p^{m/2}+1$
and a divisor  $r$ of $m$ with $1\leq r \leq m/2$ there exists a maximal
curve $C^{\prime}$ over
$\F_{p^m}$ with
$$
\eqalign{
g(C^{\prime})& =(p^{m/2}-d+1)(p^r-1)/2d,\cr
\# C^{\prime}(\F_{p^m})&=p^m+1+(p^m-(d-1)p^{m/2})(p^r-1)/d.\cr}
$$
\par
\bigskip
\centerline{\bf References}
\smallskip 
\noindent
[F-T] R.\ Fuhrmann, F.\ Torres: The genus of curves over finite
fields with many rational points. {\sl Manuscripta  Math.\ \bf 89}
(1996), p.\ 103-106.
\smallskip
\noindent
[G-V1] G.\ van der Geer, M.\ van der Vlugt : Fibre
products of Artin-Schreier curves and generalized
Hamming weights of codes. {\sl J.
Comb. Theory  A \bf 70}(1995), p.\ 337-348.
\smallskip \noindent
[G-V2]  G.\ van der Geer, M.\ van der Vlugt : How to construct curves over
finite fields with many points.  In: {\sl Arithmetic
Geometry}, (Cortona 1994), F. Catanese Ed., Cambridge Univ.\ Press,
Cambridge, 1997, p.\ 169-189.
\smallskip \noindent
[G-V3]  G.\ van der Geer, M.\ van der Vlugt: Tables for the
function $N_q(g)$. Regularly updated tables  available at {\hbox { { \tt
http://www.wins.uva.nl/ $ \tilde{ }$geer }} }. Version August 1997.
\smallskip
\noindent
[H-P] P.\ Heijnen, R.\ Pellikaan: Generalized Hamming weights of
$q$-ary Reed-Muller codes. Report Eindhoven  University of
Technology, Eindhoven 1996.
\smallskip
\noindent
[L] G.\ Lachaud: Sommes d'Eisenstein et nombre de points de
certaines courbes alg\'e\-bri\-ques sur les corps finis. {\sl C.R.\ Acad.\ Sci.\
Paris \bf 305}, S\'erie I (1987), p.\ 729-732.
\smallskip
\noindent 
[S-X] H.\ Stichtenoth, C. Xing: The genus of maximal function
fields. {\sl Manuscripta Math. \bf 86} (1995), p.\ 217-224.
\smallskip
\noindent 
[W] V.K.\ Wei: Generalized Hamming weights for linear codes. {\sl IEEE Trans.\
Inform.\ Th.\ \bf 37} (1991), p.\ 1412-1418.
\bigskip
\settabs3 \columns
\+G. van der Geer  &&M. van der Vlugt\cr
\+Faculteit
Wiskunde en Informatica &&Mathematisch Instituut\cr
\+Universiteit van
Amsterdam &&Rijksuniversiteit te Leiden \cr
\+Plantage Muidergracht 24&&Niels Bohrweg 1 \cr
\+1018 TV Amsterdam
&&2300 RA Leiden \cr
\+The Netherlands &&The Netherlands \cr
\+{\tt geer@wins.uva.nl} &&{\tt vlugt@wi.leidenuniv.nl} \cr
\bye